\documentclass[10pt]{article}
\usepackage{epsfig,makeidx,amsmath,amsfonts,latexsym,subfigure}
\newtheorem{theorem}{Theorem}[section]
\newtheorem{prop}{Proposition}[section]
\newtheorem{lemma}{Lemma}[section]
\newtheorem{corr}{Corollary}[section]
\newtheorem{defi}{Definition}[section]
\newcommand{\E}{{\mathbb E}}

\begin{document}

\title{Stationary random graphs on $\mathbb{Z}$ with prescribed iid
degrees and finite mean connections}

\author{Maria Deijfen \thanks{Stockholm University.
E-mail: mia@math.su.se} \and Johan Jonasson
\thanks{Chalmers University of Technology. E-mail:
jonasson@math.chalmers.se}}

\date{October 2006}

\maketitle

\thispagestyle{empty}

\begin{abstract}

\noindent Let $F$ be a probability distribution with support on
the non-negative integers. A model is proposed for generating
stationary simple graphs on $\mathbb{Z}$ with degree distribution
$F$ and it is shown for this model that the expected total length
of all edges at a given vertex is finite if $F$ has finite second
moment. It is not hard to see that any stationary model for
generating simple graphs on $\mathbb{Z}$ will give infinite mean
for the total edge length per vertex if $F$ does not have finite
second moment. Hence, finite second moment of $F$ is a necessary
and sufficient condition for the existence of a model with finite
mean total edge length.

\vspace{0.5cm}


\vspace{0.5cm}

\end{abstract}

\section{Introduction}

\noindent In the simplest random graph model, given a set of $n$
vertices, an edge is drawn independently between each pair of
vertices with some probability $p$. This model goes back to
Erd\H{o}s and R\'{e}nyi (1959) and dominated the field of random
graphs for decades after its introduction. However, during the
last few years there has been a growing interest in the use of
random graphs as models for various types of complex network
structures, see e.g.\ Newman (2003) and the references therein. In
this context it has become clear that the Erd\H{o}s-R\'{e}nyi
graph fails to reflect a number of important features of real-life
networks. For instance, an important quantity in a random graph is
the degree distribution, and, in an $n$-vertex Erd\H{o}s-Renyi
graph, if the edge probability is scaled by $1/n$, the vertex
degree is asymptotically Poisson distributed. In many real-life
networks however, the degree sequence has been observed to follow
a power law, that is, the number of vertices with degree $k$ is
proportional to $k^{-\tau}$ for some exponent $\tau>1$. Hence, the
Erd\H{o}s-Renyi graph provides a bad model of reality at this
point.\medskip

\noindent The fact that the Erd\H{o}s-Renyi model cannot give
other degree distributions than Poisson has inspired a number of
new graph models that take as input a probability distribution $F$
with support on the non-negative integers and give as output a
graph with degree distribution $F$. The most well-known model in
this context is the so called configuration model, formulated in
Wormald (1978) and later studied by Molloy and Reed (1995,1998)
and van der Hofstad et.\ al.\ (2005) among many others. It works
so that each vertex is assigned a random number of stubs according
to the desired degree distribution $F$ and these stubs are then
paired at random to form edges. A drawback with the configuration
model is that it can give self-loops and multiple edges between
vertices, something that in most cases is not desirable in the
applications. See however Britton et.\ al.\ (2005) for
modifications of the model that give simple graphs -- that is,
graphs without self-loops and multiple edges -- as final result. A
different model for generating graphs with given degrees is
described in Chung and Lu (2002:1,2).\medskip

\noindent A natural generalization of the problem of generating a
random graph with a prescribed degree distribution is to consider
spatial versions of the same problem, that is, to ask for a way to
generate edges between vertices arranged on some spatial structure
so that the vertex degrees have a certain specified distribution.
This problem was introduced in Deijfen and Meester (2006), where
also a model is formulated for generating stationary graphs on
$\mathbb{Z}$ with a given degree distribution $F$. Roughly the
model works so that a random number of ``stubs'' with distribution
$F$ is attached to each vertex. Each stub is then randomly and
indepedently of other stubs assigned a direction, left or right,
and the graph is obtained by stepwise pairing stubs that point to
each other. Under the assumption that $F$ has finite mean, this is
shown to lead to well-defined configurations, but the expected
length of the edges is infinite. It is conjectured that, in fact,
all stationary procedures for pairing stubs with random,
independent directions give connections with infinite mean. This
conjecture has been proved for stub distributions $F$ with bounded
support.\medskip

\noindent The purpose of the present paper is to formulate a model
for generating stationary simple graphs on $\mathbb{Z}$ with
prescribed degree distribution and finite expected edge length.
Just as in Deijfen and Meester (2006) we will begin by attaching
stubs to the vertices according to the desired degree distribution
and then we will look for a stationary way to pair these stubs to
create edges. The aim is to do this in such a way that multiple
edges are avoided. Furthermore, in view of the conjecture from
Deijfen and Meester, if we want to achieve finite mean for the
edge length, the pairing step cannot involve giving independent
directions to the stubs, but the stubs have to be connected in a
more ``effective'' way. For distributions with bounded support it
turns out that there is a quite simple way of doing this, while
the case with unbounded support requires a bit more work. We
mention that related matching problems have been studied for
instance by Holroyd and Peres (2003,2005).\medskip

\noindent Let $D\sim F$ be the degree of the origin in a
stationary simple graph on $\mathbb{Z}$ and write $T$ for the
total length of all edges at the origin. We then have that $T\geq
2\sum_{k=1}^{\lfloor D/2 \rfloor}k$, where the lower bound is
attained when there is one edge to each nearest neighbor, one edge
to each second nearest neighbor, and so on. Since
$2\sum_{k=1}^{\lfloor D/2 \rfloor}k\geq (D^2-1)/4$, it follows
that finite second moment of the degree distribution is a
necessary condition for the possibility of generating a stationary
simple graph with finite mean total edge length per vertex. In
this paper we propose a model that indeed gives finite mean for
the total edge length when $\E[D^2]<\infty$. This establishes the
following theorem.

\begin{theorem}
Let $F$ be a probability distribution with support on the
non-negative integers. It is possible to generate a simple
stationary graph on $\mathbb{Z}$ with degree distribution $F$ and
$\E[T]<\infty$ if and only if $F$ has finite second moment.
\end{theorem}

\noindent The rest of the paper is organized as follows. In
Section 2, a model is described that works for degree
distributions with bounded support, that is, for distributions
with bounded support the expected total length of all edges at a
given vertex is finite. In Section 3, this model is refined in
that vertices with high degree are treated separately. Finally in
Section 4, the refined model is shown to give finite mean for the
total edge length per vertex if the degree distribution has finite
second moment.

\section{A basic model}

\noindent Let $F$ be a probability distribution with support on
the non-negative integers. In this section we formulate a basic
model for generating a stationary simple graph on $\mathbb{Z}$
with degree distribution $F$. We also show that the expected
length of the edges is finite if $F$ has bounded support.\medskip

\noindent The basis for the model proposed in this section -- and
also for the refined model in the following section -- is a stub
configuration on $\mathbb{Z}$. This configuration is obtained by
associating independently to each vertex $i\in\mathbb{Z}$ a random
degree $D_i$ with distribution $F$ and then attach $D_i$ stubs to
vertex $i$. The question is how the stubs should be connected to
create edges. As mentioned, an important restriction is that the
pairing procedure is required to be stationary. Also, the
resulting graph is not allowed to contain multiple edges between
vertices.\medskip

\noindent Our first suggestion of how to join the stubs is as
follows. Let $\Gamma_j$ be the set of all vertices with degree at
least $j$, that is, $\Gamma_j=\{i\in\mathbb{Z}: D_i\geq j\}$.
Furthermore, for each vertex $i$, label the stubs
$\{s_{i,j}\}_{j=1}^{D_i}$ and define $\Lambda_j=\bigcup_{i \in
\Gamma_j}s_{i,j}$ so that each vertex $i\in\Gamma_j$ has exactly
one stub $s_{i,j}$ in $\Lambda_j$ connected to it. Stubs in
$\Lambda_j$ will be referred to as belonging to \emph{level} $j$.
The stubs in the sets $\{\Lambda_j\}$ are now connected to each
other within the sets, starting with $\Lambda_1$, as follows:

\begin{itemize}
\item[1.] Imagine that a fair coin is tossed at the first vertex
$i_1\geq 0$ with degree at least 1. If the coin comes up heads,
the level 1 stub $s_{i_1,1}$ at $i_1$ is turned to the right and
if the coin comes up tails, the stub $s_{i_1,1}$ is turned to the
left. The other stubs in $\Lambda_1$ are directed so that, at
every second vertex in $\Gamma_1$, the level 1 stub points to the
right and, at every second vertex, it points to the left. Edges
are then created by connecting stubs that are pointing at each
other. More precisely, a level 1 stub at vertex $i$ pointing to
the right (left) is joined to the level 1 stub at the next vertex
$j>i$ ($j<i$) in $\Gamma_1$ -- by the construction, this stub is
pointing to the left (right).

\item[2.] The stubs in $\Lambda_2$ are assigned directions
analogously by letting every second stub in $\Lambda_2$ point to
the right and every second stub to the right, the direction of the
stub at the first vertex $i_2\geq 0$ with degree at least 2 being
determined by a coin toss. To avoid multiple edges between
vertices, the stubs are then connected so that a right (left) stub
$s_{i,2}$ at a vertex $i$ is joined with the \emph{second} left
(right) stub encountered to the right (left) of $i$. Since every
second stub in $\Lambda_2$ is pointing to the right (left), this
means that $s_{i,2}$ is linked to the level 2 stub at the third
vertex in $\Gamma_2$ to the right of $i$. This vertex cannot have
an edge to $i$ from step 1, since the level 1 stub at $i$ was
connected to the first vertex in $\Gamma_1$ either to the left or
to the right, and $\Gamma_2\subset\Gamma_1$.
$$
\vdots
$$
\item[$n$.] In general, in step $n$, the stubs in $\Lambda_n$ are
connected by first randomly choosing one of the two possible
configurations where every second stub is pointing to the right
and every second stub is pointing to the left, and then link a
given stub to the $n$th stub pointing in the opposite direction
encountered in the direction of the stub. A level $n$ stub at
vertex $i$ pointing to the right (left) is hence connected to
vertex number $2n-1$ to the right (left) of $i$ in $\Gamma_n$.
Since $\Gamma_n\subset\Gamma_{n-1}$, this does not give rise to
multiple edges.
$$
\vdots
$$
\end{itemize}

\noindent This procedure is clearly stationary and will be
referred to as the \emph{Coin Toss} (CT) \emph{model}. Our first
result is a formula for the expected total edge length per vertex
in the CT-model. To formulate it, write $p_j$ for the probability
of the outcome $j$ in the degree distribution $F$ and note that,
by stationarity, the distribution of the edge length is the same
at all vertices. Hence it suffices to consider the total length
$T$ of all edges at the origin.

\begin{prop}\label{prop:msct}
In the CT-model, assume that $F$ has bounded support and let
$u=\max\{j:p_j>0\}$. Then
$$
\E[T]=u^2.
$$
\end{prop}

\noindent \textbf{Proof:} In what follows we will drop the vertex
index for the degree at the origin and write $D_0=D$. For
$j=1,\ldots,u$, let $K_j$ denote the length of the $j$th edge at
the origin, that is, $K_j$ is the length of the edge created by
the stub $s_{0,j}$ (if $D<j$, we set $K_j=0$). Also, define
$p_j^+=P(D\geq j)=1-F(j-1)$. Then

\begin{eqnarray*}
\E[T] & = & \E\Big[\sum_{j=1}^uK_j\Big]\\
& = & \sum_{j=1}^up_j^+ \E[K_j|D \geq j].
\end{eqnarray*}

\noindent The $j$th edge at the origin is equally likely to point
to the right as to the left, and, by symmetry, the expected length
of the edge is the same in both cases. If the edge points to the
right (left), its other endpoint is vertex number $2j-1$ to the
right (left) of the origin with degree at least $j$. The distance
to this vertex has a negative binomial distribution with mean
$(2j-1)/p_j^+$, and it follows that
$$
\E[T]=\sum_{n=1}^u\frac{2j-1}{p_j^+}p_j^+ = 2\sum_{j=1}^uj-u =
u^2,
$$
as desired.\hfill$\Box$\medskip

\noindent It follows from the calculations in the proof that $T$
has infinite mean when the support of $F$ is unbounded and hence
we have the following corollary:

\begin{corr}\label{corr}
In the CT-model $\E[T]<\infty$ iff $F$ has bounded support.
\end{corr}

\noindent Roughly, the reason the CT-model gives infinite mean for
$T$ when the degrees are not bounded is that vertices with high
degree are connected to other vertices with high degree. In Section 3, a model is formulated where high degree vertices are connected in a more effective way in order to get edges with finite mean length.

According to Corollary \ref{corr}, the mean total edge length at
the origin in the CT-model with i.i.d.\ degrees is finite
for degree distributions with bounded support. The following
proposition -- which will be needed in the proof of Theorem
\ref{th:fmcm} -- asserts that the same result holds also for
stationary degrees.

\begin{prop}\label{prop:vvcts}
Let $\{D'_i\}$ be a stationary stub configuration on $\mathbb{Z}$
with $D'_i\sim G$ and connect the stubs as in the CT-model. Then,
if $G$ has bounded support, we have $\E[T]<\infty$.
\end{prop}

The proof of the proposition is based on the following lemma.

\begin{lemma}\label{lemma:statvv}
Let $\{X_i\}_{i\in\mathbb{Z}}$ be a $\{0,1\}$-valued stationary
process with an almost surely infinite number of 1's. Define
$\tau_i=\inf\{k>i:\, X_k=1\}$ and write $P(X_i=1)=p$. Then
$$
\E[\tau_i|X_i=1]=\frac{1}{p}.
$$
\end{lemma}

\noindent \textbf{Proof of Lemma \ref{lemma:statvv}:} By
stationarity, it suffices to consider $i=0$. We have
$$
\E[\tau_0|X_0=1]=\frac{\E[\tau_01_{\{X_0=1\}}]}{p},
$$
where the numerator can be written as

\begin{eqnarray*}
\E[\tau_01_{\{X_0=1\}}] & = & \sum_{k=0}^\infty
P(\tau_01_{\{X_0=1\}} > k)\\
& = & \sum_{k=0}^\infty P(X_0=1,X_1=0,\ldots ,X_k=0)\\
& = & \sum_{k=0}^\infty P(X_{-k}=1,X_{-(k-1)}=0,\ldots ,X_0=0) \\
& = & \sum_{k=0}^\infty P(\min\{j \geq 0:X_{-j}=1\} = k) \\
& = & 1.
\end{eqnarray*}
\hfill$\Box$\medskip

\noindent \textbf{Proof of Proposition \ref{prop:vvcts}:} Write
$D_0'=D'$, let $u=\max\{j:p_j>0\}$, and remember from the proof of
Proposition \ref{prop:msct} that $K_j$ denotes the length of the
$j$th edge at the origin for $j=1,\ldots,u$ (if $D'<j$, we set
$K_j=0$). Clearly we are done if we can show that E$[K_j]<\infty$
for all $j$. We have E$[K_j]=\E[K_j|D'\geq j]p_j^+$, and to see
that $\E[K_j|D'\geq j]$ is finite, define
$$
X_i^j = \left\{ \begin{array}{ll}
                      1 & \mbox{if $D'\geq j$};\\
                      0 & \mbox{otherwise}.
                    \end{array}
            \right.
$$

\noindent The $j$th edge at the origin is equally likely to point
to the right as to the left, and, by symmetry, the expected length
of the edge is the same in both cases. If the edge points to the
right (left), its other endpoint is vertex number $2j-1$ to the
right (left) of the origin with degree at least $j$. Hence
$\E[K_j|D'\geq j]$ is the expected distance to the right of the
origin until $2j-1$ 1's have been encountered in the process
$\{X_i^j\}$, given that $X_0=1$. It follows from Lemma
\ref{lemma:statvv} that the expected distance between two
successive 1's in the process is $1/p_j^+$, which gives the
desired result exactly as in the proof of Proposition
\ref{prop:msct}.\hfill$\Box$\medskip

\section{A refined model}

\noindent In this section we describe a model designed to give
finite mean for the total edge length per vertex for degree
distribution with finite second moment. The idea is to truncate
the degrees at some high level $d$ and connect stubs at level
$j\geq d+1$ separately. The remaining stubs at level $1,\ldots,d$
after this has been done are then connected according to the
CT-model. The model is slightly easier to define when the degrees
are almost surely non-zero, and hence, in what follows, we will
assume that $P(D_i\geq 1)=1$. This means no loss of generality,
since removing vertices with degree 0 only shrinks expected edge
lengths by a factor $1-P(D_i=0)$.\medskip

\noindent First we introduce some notation and terminology. For a
fixed $d\in\mathbb{N}$, let $D_i^{d}$ be the ``tail'' above level
$d$ at vertex $i$, that is, $D_i^d =\max\{D_i-d,0\}$. Stubs at
level $j\geq d+1$ will be referred to as \emph{bad} and $D_i^{d}$
thus indicates the number of bad stubs at vertex $i$. A vertex
with bad stubs on it -- that is, a vertex with degree strictly
larger than $d$ -- will be called \emph{high}. Of course, stubs
that are not bad will be called {\em good} and vertices that are
not high will be called {\em low}. The model for connecting the
stubs is based on the concept of \emph{claimed} vertices:

\begin{defi}
A vertex $i\in\mathbb{Z}$ is said to be claimed on level $d$ iff
$\sum_{k=i-m}^{i+m}D_k^{d}\geq m$ for some $m\geq 1$.
\end{defi}

\noindent In words, a vertex is claimed if, either there are bad
stubs on the vertex itself (this means that a high vertex is by
definition claimed), or there is a symmetric interval of width
$2m+1$ for some positive $m$ around the vertex in which the total
number of bad stubs is at least $m$. The set of claimed vertices
at level $d$ will be denoted by $\mathcal{C}_d$. Finally, by a
\emph{cluster} of claimed vertices we mean a set of consecutive
vertices $i,\ldots,i+n$ with $\{i,\ldots,i+n\}\subset
\mathcal{C}_d$ but $i-1\not\in \mathcal{C}_d$ and $i+n+1\not\in
\mathcal{C}_d$.
\medskip

\noindent The idea with these definitions is that a claimed vertex
that is not high is in some sense close to a high vertex and might
therefore be used for the bad stubs at the high vertices to
connect to. Indeed, in the model that we will soon propose for
connecting the stubs, bad stubs are always connected within the
claimed cluster that their vertex belongs to. To make sure that
this can be done without creating multiple edges we need to see
that the total number of bad stubs in a claimed cluster is
strictly smaller than the number of vertices in the cluster.
Hence, for a given subset $A$ of $\mathbb{Z}$, let $b(A)$ be the
number of bad stubs at vertices in $A$, that is,
$$
b(A)=\sum_{i\in A}D_i^{d}.
$$
Then the following holds:

\begin{lemma}\label{lemma:tmdsik}
For each cluster $C$ in $\mathcal{C}_d$, we have $b(C)\leq |C|-1$.
\end{lemma}

\noindent \textbf{Proof:} Consider a given claimed cluster
$C=[i+1,i+n]$ and assume that $b(C)\geq |C|$, that is, $b(C)\geq
n$. We then have for the vertex $i$ next to the left endpoint of
the cluster that
$$
\sum_{k=i-n}^{k=i+n}D_i^{d}\geq b(C)\geq n
$$
so that hence $i$ is claimed as well, which is a
contradiction.\hfill$\Box$\medskip

\noindent We are now ready to describe the refined model. The model
requires that the claimed clusters are almost surely finite. This
will indeed be the case if $d$ is large, as follows from
Proposition \ref{prop:eca} in the next section (which stipulates
that the expected cluster size is finite for large $d$). For now
we take this for granted. Hence, fix $d$ large enough to ensure
that the clusters are finite and, to make step 2 below easier,
assume without loss of generality that $d$ is even.

\begin{itemize}
\item[1.] Bad stubs are connected within the claimed clusters. In
a cluster with only one single high vertex $i$, this is done by
choosing a random subset of $D_i^{d}$ low vertices in the cluster,
pick one stub from each of these vertices and connect to a bad
stub of $i$. This is indeed possible, since, by Lemma
\ref{lemma:tmdsik}, there are at least $D_i^{d}$ low vertices in
the cluster (in fact, at least $2D_i^d$) and, by assumption, there
is at least one stub at each vertex. If there is more than one
high vertex in a cluster $C$, the bad stubs are connected as
follows.

\begin{itemize}
\item[(i)] Write $h=h(C)$ for the number of high vertices in $C$
and let these high vertices be denoted by $i_1,\ldots,i_{h}$
(ordered from the left to the right). First we use some of the bad
stubs to create edges between high vertices: If $h$ is even,
consider the pairs $(i_1,i_2),\ldots,(i_{h-1},i_{h})$ and connect
the two vertices in each pair by using one bad stub from each
vertex. If $h$ is odd, with probability 1/2, leave the last high
vertex $i_h$ out and connect the pairs $(i_1,i_2),\ldots
(i_{h-2},i_{h-1})$, and, with probability 1/2, leave the first
high vertex $i_1$ out and connect the pairs $(i_2,i_3),\ldots
(i_{h-1},i_h)$. This means that one bad stub from each high vertex
(except $i_h$ or $i_1$) is connected if $h$ is even (odd). In any
case, at least $h-1$ bad stubs are used.

\item[(ii)] Let $b_r(C)$ denote the number of remaining
unconnected bad stubs in $C$ after step (i). By Lemma
\ref{lemma:tmdsik}, the total number of bad stubs in $C$ is at
most $|C|-1$, implying that $b_r(C)\leq |C|-h$, that is, the
number of unconnected bad stubs in $C$ does not exceed the number
$|C|-h$ of low vertices. Hence, to connect the remaining bad
stubs, chose randomly $b_r(C)$ low vertices, take one stub from
each of these vertices and pair these stubs randomly with the bad
stubs.
\end{itemize}

\item[2.] In this step, the good stubs at the high vertices
$\{i_j\}$ are connected to each other. Each high vertex has $d$
good stubs attached to it and, for a given high vertex $i_j$,
these stubs are linked to good stubs at the vertices
$i_{j+2},\ldots,i_{j+d/2+1}$ and $i_{j-2},\ldots,i_{j-d/2-1}$,
that is, half of the good stubs are pointed to the left and half
of them to the right. Consecutive high vertices might already have
an edge between them from step 1(i) and therefore, to avoid
multiple edges, we do not connect $i_j$ to $i_{j-1}$ or $i_{j+1}$.

\item[3.] The remaining stubs at low vertices are connected
according to the CT-model. There are no edges between low vertices
from the previous steps and hence multiple edges will not arise
here.
\end{itemize}

\noindent This model will be referred to as the \emph{cluster
model}. The next task is to show that the mean total edge length
per vertex is finite provided the degree distribution has finite
second moment.

\section{The mean total edge length per vertex}

\noindent The truncation level $d$ is of course important for the
properties of the cluster model. Write $T^{d}$ for the total
length of all edges at the origin for a given value of $d$. The
aim in this section is to prove the following theorem.

\begin{theorem}\label{th:fmcm}
If $F$ has finite second moment, then, for large $d$, we have
$\E[T^{d}]<\infty$ in the cluster model.
\end{theorem}

\noindent A large part of the work in proving this theorem lies in
showing that the expected size of a claimed clusters is finite if
$d$ is large. Since the bad stubs are connected within the claimed
clusters, this ensures that the expected length of the edges
created by the bad stubs is finite. For technical reasons, we will
need a slightly more general result. To formulate it, first
generalize the definition of a claimed vertex to incorporate a
parameter $\alpha$.

\begin{defi}
A vertex is $\alpha$-claimed on level $d$ iff
$\sum_{k=i-m}^{i+m}D_k^{d}\geq\alpha m$ for some $m\geq 1$.
\end{defi}

\noindent Given a stub configuration $\{D_i\}_{i\in\mathbb{Z}}$,
we can now talk about clusters of $\alpha$-claimed vertices. Let
$C^{d,\alpha}$ be the $\alpha$-claimed cluster of the origin.
Also, let $\mu_d=\E[D_i^d]$, that is, $\mu_d$ is the expected
number of stubs above level $d$ at a given vertex. Clearly
$\mu_d\rightarrow 0$ as $d\rightarrow\infty$ and hence, by picking
$d$ large, we can make $\mu_d$ arbitrarily small. The result
concerning the expected cluster size now runs as follows.

\begin{prop}\label{prop:eca}
Fix $\alpha>0$. If $d$ is large enough to ensure that
$\mu_d<\alpha/18$, then $\E\left[|C^{d,\alpha}|\right]<\infty$.
\end{prop}

\noindent\textbf{Proof:} We will show that

\begin{equation}\label{eq:kosn}
\big\{|C^{d,\alpha}|\geq n\big\}\subset\Big\{\exists m\geq
n:\,\sum_{k=-m}^m D_k^d\geq \frac{\alpha m}{6}\Big\}.
\end{equation}

\noindent With $\widetilde{D}_k^d=D_k^d-\mu_d$, this implies that

\begin{eqnarray*}
\big\{|C^{d,\alpha}|\geq n\big\} &\subset & \Big\{\exists m\geq
n:\,\sum_{k=-m}^m \widetilde{D}_k^d\geq \frac{\alpha
m}{6}-(2m+1)\mu_d\Big\}\\
& \subset & \Big\{\exists m\geq n:\,\sum_{k=-m}^m
\widetilde{D}_k^d\geq c_{d,\alpha}m\Big\},
\end{eqnarray*}

\noindent where $c_{d,\alpha}:=\alpha/6-3\mu_d>0$ (to get the last
inclusion, we have used that $2m+1\leq 3m$). Hence

\begin{eqnarray*}
P\big(|C^{d,\alpha}|\geq n\big) &\leq & P\left(\exists m\geq
n:\,\sum_{k=-m}^m \widetilde{D}_k^d\geq c_{d,\alpha}m\right)\\
& = & P\left(\sup_{m\geq n}\,\frac{1}{m}\sum_{k=-m}^m
\widetilde{D}_k^d\geq c_{d,\alpha}\right),
\end{eqnarray*}

\noindent and consequently
$$
\E\big[|C^{d,\alpha}|\big]\leq \sum_{n=1}^\infty
P\left(\sup_{m\geq n}\,\Big(\frac1m\sum_{k=-m}^m
\widetilde{D}_k^d\Big)\geq c_{d,\alpha}\right).
$$
By a standard result on convergence rate in the law of large
numbers from Baum and Katz (1965; Theorem 3 with $t=r=2$), the sum on the
right hand side is convergent iff the $\widetilde{D}_k^d$'s have
finite second moment. This proves the proposition.\medskip

\noindent It remains to show (\ref{eq:kosn}). To this end, for
$i\in C^{d,\alpha}$, let $\mathcal{I}_i$ be the shortest interval
around $i$ where the condition for $i$ to be $\alpha$-claimed is
satisfied. More precisely, if
$$
m_i=\inf\Big\{m:\, \sum_{k=i-m}^{i+m}D_k^d\geq \alpha m\Big\},
$$
we have $\mathcal{I}_i=[i-m_i,i+m_i]$. We now claim that we can
pick a subset $\{\mathcal{I}_{i_j}\}_{i_j\in C^{d,\alpha}}$ of
these intervals that completely covers the cluster
($C^{d,\alpha}\subset \cup_j\mathcal{I}_{i_j}$) and where only
consecutive intervals intersect ($\mathcal{I}_{i_j}\cap
\mathcal{I}_{i_k}=\emptyset$ if $|j-k|\geq 2$). To construct this
subsequence, let $\mathcal{I}_{i_1}$ be the interval in
$\{\mathcal{I}_i\}_{i\in C^{d,\alpha}}$ that reaches furthest to
the left, that is, $\mathcal{I}_{i_1}$ is the interval with
$l=\inf\{k:\,k\in\cup_{i\in C^{d,\alpha}}\mathcal{I}_i\}$ as its
left endpoint. If there is more than one interval in
$\{\mathcal{I}_i\}_{i\in C^{d,\alpha}}$ with $l$ as its left
endpoint, we take $\mathcal{I}_{i_1}$ to be the largest one. Next,
consider the set $S_1$ of intervals $\mathcal{I}_k$ with $k\in
C^{d,\alpha}$ that intersect $\mathcal{I}_{i_1}$ and define
$\mathcal{I}_{i_2}$ to be the interval in $S_1$ that reaches
furthest to the right. If there is more than one interval in $S_1$
that ends at the same maximal right endpoint, we let
$\mathcal{I}_{i_2}$ be the largest of those intervals. Let $S_2$
be the set of intervals that intersect $\mathcal{I}_{i_2}$. The
interval $\mathcal{I}_{i_3}$ is set to be the member in $S_2$ that
reaches furthest to the right and, as before, if there is more
than one candidate, we pick the largest one. This interval cannot
intersect $\mathcal{I}_{i_1}$, since then it would have been
chosen already in the previous step when $\mathcal{I}_{i_2}$ was
defined. In general, given $\mathcal{I}_{i_1},\ldots,
\mathcal{I}_{i_j}$ the interval $\mathcal{I}_{i_{j+1}}$ is defined
as follows.

\begin{itemize}
\item[(i)] Let $S_j=\{\mathcal{I}_{i}:\,i\in C^{d,\alpha}\textrm{
and }\mathcal{I}_{i}\cap\mathcal{I}_{i_j}\neq\emptyset\}$ and
write $r_j=\sup\{k:\,k\in S_j\}$.

\item[(ii)] Take $\mathcal{I}_{i_{j+1}}$ to be the largest
interval in $S_j$ with its right endpoint at the vertex $r_j$.
\end{itemize}

\noindent We repeat this procedure until an interval
$\mathcal{I}_{i_s}$ whose right endpoint is outside the cluster
$C^{d,\alpha}$ is picked. The entire cluster is then covered by
$\cup_{j=1}^s\mathcal{I}_{i_j}$ and, by construction,
non-consecutive intervals do not intersect, as desired.\medskip

\noindent Now let $A=\cup_{k\geq 1}\mathcal{I}_{i_{2k-1}}$ and
$B=\cup_{k\geq 1}\mathcal{I}_{i_{2k}}$, that is, every second
interval in $\{\mathcal{I}_{i_j}\}$ is placed in $A$ and every
second interval is placed in $B$. Then $A$ and $B$ are both unions
of mutually disjoint intervals. Remember that $b(\cdot)$ denotes
the number of bad stubs in a given set and note that, by the
definition of badness, for a given interval $\mathcal{I}_{i_j}$,
we have $b(\mathcal{I}_{i_j})\geq \alpha|\mathcal{I}_{i_j}|/3$.
Since the intervals in $A$ are disjoint, it follows that $b(A)\geq
\alpha|A|/3$, and similarly, $b(B)\geq \alpha|B|/3$. Hence

\begin{eqnarray*}
b(A\cup B) & \geq & \max\{b(A),b(B)\}\\
& \geq & \frac{\alpha}{3}\max\{|A|,|B|\}\\
& \geq & \frac{\alpha}{6}|A\cup B|.
\end{eqnarray*}

\noindent With $m=|A\cup B|$, we have

\begin{eqnarray*}
b([-m,m]) & \geq & b(A\cup B)\\
& \geq & \frac{\alpha}{6}|A\cup B|\\
& = & \frac{\alpha m}{6}.
\end{eqnarray*}

\noindent This establishes (\ref{eq:kosn}).\hfill$\Box$\medskip

\noindent We are now ready to prove Theorem \ref{th:fmcm}. The
proof is based on Proposition \ref{prop:eca}, which ensures that
the expected length of the edges created by the bad stubs in step
1 in the description of the cluster model is finite, and
\ref{prop:vvcts}, which guarantees that the edges created in step
2 and 3 have finite mean.\medskip

\noindent \textbf{Proof of Theorem \ref{th:fmcm}:} Write $T_1^d$,
$T_2^d$ and $T_3^d$ for the total length of the edges created at
the origin in step 1, 2 and 3 respectively in the description of
the cluster model.\medskip

\noindent First we attack $T_1^d$. To this end, given a stub
configuration $\{D_i\}_{i\in\mathbb{Z}}$, write $C^{d,\alpha}_i$
for the $\alpha$-claimed cluster of the vertex $i$. Now, $T_1^d$
is the total length of all edges created by bad stubs at the
origin. The number of such edges is clearly smaller than the total
number $D$ of edges at the origin and they are all connected
within the claimed cluster of the origin. Hence $T_1^d\leq
D|C^{d,1}_0|$ (the cluster model is based on $\alpha=1$). If $D$
and $|C^{d,1}_0|$ were independent it would follow immediately
from Proposition \ref{prop:eca} that E$[T_1^d]<\infty$ for large
$d$. However, $D$ and $|C^{d,1}_0|$ are of course not independent.
To get around this, introduce a coupled degree configuration
$\{\widehat{D}_i\}_{i\in\mathbb{Z}}$ where $\widehat{D}$ is
generated independently, while $\widehat{D}_i=D_i$ for $i\neq 0$.
Quantities based on $\{\widehat{D}_i\}_{i\in\mathbb{Z}}$ will be
equipped with a hat-symbol. We will show that

\begin{equation}\label{eq:ksb}
\big|C^{d,1}_0\big|\leq 4D+\big|\widehat{C}^{d,1/2}_{-2D}\big|
+\big|\widehat{C}^{d,1/2}_{2D}\big|.
\end{equation}

\noindent Since $\widehat{C}^{d,1/2}_i$ is clearly independent of
$D$ for all $i$, this implies that

\begin{eqnarray*}
\E[T^d_1] & \leq &
\E\left[D\left(4D+\big|\widehat{C}^{d,1/2}_{-2D}\big|
+\big|\widehat{C}^{d,1/2}_{2D}\big|\right)\right]\\
& = &4\E\left[D^2\right]+2\E\left[D\right]
\cdot\E\left[\big|\widehat{C}^{d,1/2}_{2D}\big|\right].
\end{eqnarray*}

\noindent If $F$ has finite second moment and $d$ is large so that
$\mu_d\leq 1/36$, then, by Proposition \ref{prop:eca}, we have
E$\big[\big|\widehat{C}^{d,1/2}_{2D}\big|\big]<\infty$. It follows
that E$[T^d_1]$ is finite under the same conditions.\medskip

\noindent To establish (\ref{eq:ksb}), it suffices to observe that
each vertex $i\not\in [-2D-1,2D+1]$ that is claimed for $\alpha=1$
in the original configuration $\{D_i\}$ is still claimed for
$\alpha=1/2$ in the coupled configuration $\{\widehat{D}_i\}$.
Hence, pick a vertex $i$ with $|i|\geq 2D$ that is claimed for
$\alpha=1$ in the original configuration. Write
$\mathcal{I}_i^m=[i-m,i+m]$ for the smallest interval such that
$b(\mathcal{I}_i^m)\geq m$ and assume that $m\geq 2D$ so that
hence $0\in \mathcal{I}_m$ (if this is not the case, $i$ is
obviously claimed for $\alpha=1$ also in the coupled configuration
$\{\widehat{D}_i\}$, since $\widehat{D}_i=D_i$ for all $i\neq 0$).
For such $i$, we have

\begin{eqnarray*}
\widehat{b}(\mathcal{I}_i^m) & \geq &
b(\mathcal{I}_i^m)-D\\
&\geq & m-D\\
& \geq & m/2,
\end{eqnarray*}

\noindent meaning that $i$ is claimed for $\alpha=1/2$ in
$\{\widehat{D}_i\}$, as desired.\medskip

\noindent Next, consider the total length $T^d_2$ of the edges
created at the origin in step 2, where good stubs at high vertices
are connected. If the origin is not high, that is, if $D\leq d$,
then clearly $T^d_2=0$. Hence assume that $D\geq d+1$. Then $d$
edges will be created at the origin in step 2 -- half of them will
point to the right and half of them to the left. The longest edge
to the right (left) runs to vertex number $2d+1$ to the right
(left) of the origin with degree larger than or equal to $d+1$.
The distance to this vertex has a negative binomial distribution
with finite mean, and it follows that $T^d_2$ has finite
mean.\medskip

\noindent All that remains is to see that the total length $T^d_3$
of the edges created in step 3 -- when the CT-model is applied to
connect remaining stubs after steps 1 and 2 -- has finite
expectation. This however is an immediate consequence of
Proposition \ref{prop:vvcts}, since, if $D_i'$ denotes the number
of stubs at vertex $i$ that are not connected after steps 1 and 2,
then $D_i'\leq d$ and $\{D'_i\}$ is clearly a stationary
sequence.\medskip

\noindent To sum up, we have shown that $\E[T]=
\E[T_1^d+T_2^d+T_3^d]<\infty$, as desired.\hfill$\Box$\bigskip

\noindent \textbf{Acknowledgement} We thank Olle H\"{a}ggstr\"{o}m
for giving the idea for the coin toss model.

\section*{References}

\noindent Baum, E. and Katz, M (1965): Convergence rates in the
law of large numbers, \emph{Trans. Amer. Math. Soc.} \textbf{120},
108-123. \medskip

\noindent Britton, T., Deijfen, M. and Martin-L\"{o}f, A. (2005):
Generating simple random graphs with prescribed degree
distribution, \emph{J. Stat. Phys.}, to appear.\medskip

\noindent Chung, F. and Lu, L. (2002:1): Connected components in
random graphs with given degrees sequences, \emph{Ann. Comb.}
\textbf{6}, 125-145.\medskip

\noindent Chung, F. and Lu, L. (2002:2): The average distances in
random graphs with given expected degrees, \emph{Proc. Natl. Acad.
Sci.} \textbf{99}, 15879-15882.\medskip

\noindent Deijfen, M and Meester, R. (2006): Generating stationary
random graphs on $\mathbb{Z}$ with prescribed i.i.d.\ degrees,
\emph{Adv. Appl. Probab.} \textbf{38}, 287-298.\medskip

\noindent Erd\H{o}s, P. and R\'{e}nyi, A. (1959): On random
graphs, \emph{Publ. Math.} \textbf{6}, 290-297.\medskip

\noindent Hofstad, R. van der, Hooghiemstra, G. and Znamenski, D.
(2005): Random graphs with arbitrary i.i.d. degrees, preprint
(www.win.tue.nl/$\sim$rhofstad).\medskip

\noindent Holroyd, A.E. and Peres, Y. (2003): Trees and Matchings
from Point Processes, \emph{Electr. Commun. Probab.} \textbf{8:3},
17-27.\medskip

\noindent Holroyd, A.E. and Peres, Y. (2005): Extra heads and
invariant allocations. \emph{Ann. Probab.} {\bf 33},
31-52.\medskip

\noindent Molloy, M. and Reed, B. (1995): A critical point for
random graphs with a given degree sequence, \emph{Rand. Struct.
Alg.} \textbf{6}, 161-179.\medskip

\noindent  Molloy, M. and Reed, B. (1998): The size of the giant
component of a random graphs with a given degree sequence,
\emph{Comb. Probab. Comput.}\ \textbf{7}, 295-305.\medskip

\noindent Newman (2003): The structure and function of complex
networks, \emph{SIAM Rev.} \textbf{45}, 167-256.\medskip

\noindent Wormald, N.C. (1978): Some problems in the enumeration
of labelled graphs, Doctoral thesis, Newcastle University.\medskip

\end{document}